\documentclass{article}
\usepackage[utf8]{inputenc}
\usepackage{amsmath,amsthm,amssymb,amsfonts}
\usepackage{verbatim}
\usepackage{graphicx}
\usepackage{stmaryrd} 
\usepackage{hyperref}

\usepackage[colorinlistoftodos]{todonotes}

\theoremstyle{plain}
\newtheorem{thm}{Theorem}
\newtheorem{lem}[thm]{Lemma}
\newtheorem{prop}[thm]{Proposition}
\newtheorem{cor}[thm]{Corollary}

\newtheorem{remk}[thm]{Remark}

\theoremstyle{definition}
\newtheorem{definition}[thm]{Definition}
\newtheorem{defn}[thm]{Definition}
\newtheorem{exl}[thm]{Example}

\numberwithin{thm}{section}

\newcommand{\adj}{\leftrightarrow}
\newcommand{\adjeq}{\leftrightarroweq}

\DeclareMathOperator{\id}{id}
\DeclareMathOperator{\Fix}{Fix}
\def\Z{{\mathbb Z}}
\def\N{{\mathbb N}}
\def\R{{\mathbb R}}

\title{Some Consequences of Restrictions on
Digitally Continuous Functions}
\author{Laurence Boxer
\thanks{Department of Computer and Information Sciences, Niagara University,
        New York, 14109, USA;
        and Department of Computer Science and Engineering, State University
        of New York at Buffalo \newline
        boxer@niagara.edu \newline
        To appear in {\em Note di Matematica}}
}
\date{}

\begin{document}

\maketitle
\begin{abstract}
We study the consequences of some restrictions on 
digitally continuous functions.
One of our results modifies easily to yield an 
analogous result for topological spaces.
    

\end{abstract}

\section{Introduction}
If $f: X \to Y$ is a continuous function between topological
spaces, and $\emptyset \neq A \subset X$, it is often true
that knowledge of $f|_A$ tells us little about 
$f|_{X \setminus A}$. A digital image is often a model of
an object in Euclidean space, and the concept of a digitally
continuous function is modeled on the ``preservation of nearness"
notion of a Euclidean continuous function; however, when we 
consider a continuous function 
$f: (X,\kappa) \to (Y,\lambda)$ between digital images, 
we often find that knowledge of $f|_A$ tells us much about 
$f|_{X \setminus A}$. In this paper, we continue the work
of fixed point theory for digital images
(see~\cite{Rosenfeld,BEKLL,bs19a,BxFpSets,BxConvex,BxSubsets,BxColdAndFreezing}) 
and coincidence theory for digital images
(see~\cite{AbdullahiEtAl}) by examining how restrictions
placed on $f|_A$ limit $f|_{X \setminus A}$.

\section{Preliminaries}
Let $\N$ denote the set of natural numbers; $\N^* = \{0\} \cup \N$, the
set of nonnegative integers; $\Z$, the set of integers; and $\R$, the set of real numbers.
$\#X$ will be used for the number of members of a set~$X$.

\subsection{Adjacencies}
Material in this section is largely quoted or paraphrased from~\cite{bs19a}.

A digital image is a pair $(X,\kappa)$ where
$X \subset \Z^n$ for some $n$ and $\kappa$ is
an adjacency on $X$. Thus, $(X,\kappa)$ is a graph
for which $X$ is the vertex set and $\kappa$ 
determines the edge set. Usually, $X$ is finite,
although there are papers that consider infinite $X$. Usually, adjacency reflects some type of
``closeness" in $\Z^n$ of the adjacent points.
When these ``usual" conditions are satisfied, one
may consider a subset $Y$ of $\Z^n$ (typically, an 
$n$-dimensional cube) containing $X$ as a model of a
black-and-white ``real world" image in which
the black points (foreground) are represented by 
the members of $X$ and the white points 
(background) by members of $Y \setminus \{X\}$.

We write $x \adj_{\kappa} y$, or $x \adj y$ when
$\kappa$ is understood or when it is unnecessary to
mention $\kappa$, to indicate that $x$ and
$y$ are $\kappa$-adjacent. Notations 
$x \adjeq_{\kappa} y$, or $x \adjeq y$ when
$\kappa$ is understood, indicate that $x$ and
$y$ are $\kappa$-adjacent or are equal.

The most commonly used adjacencies are the
$c_u$ adjacencies, defined as follows.
Let $X \subset \Z^n$ and let $u \in \Z$,
$1 \le u \le n$. Then for points
\[x=(x_1, \ldots, x_n) \neq (y_1,\ldots,y_n)=y\]
we have $x \adj_{c_u} y$ if and only if
\begin{itemize}
    \item for at most $u$ indices $i$ we have
          $|x_i - y_i| = 1$, and
    \item for all indices $j$, $|x_j - y_j| \neq 1$
          implies $x_j=y_j$.
\end{itemize}

The $c_u$-adjacencies are often denoted by the
number of adjacent points a point can have in the
adjacency. E.g.,
\begin{itemize}
\item in $\Z$, $c_1$-adjacency is 2-adjacency;
\item in $\Z^2$, $c_1$-adjacency is 4-adjacency and
      $c_2$-adjacency is 8-adjacency;
\item in $\Z^3$, $c_1$-adjacency is 6-adjacency,
      $c_2$-adjacency is 18-adjacency, and 
      $c_3$-adjacency is 26-adjacency.
\end{itemize}
In this paper, we mostly use the $c_1$ and $c_2$ adjacencies in $\Z^2$.

Let $x \in (X,\kappa)$. We use the notations
\[  N(X,x,\kappa) = \{ y \in X \, | \, y \adj_{\kappa} x \}
\]
and
\[  N^*(X,x,\kappa) =  \{ y \in X \, | \, y \adjeq_{\kappa} x \} = N(X,x,\kappa) \cup \{x\}.
\]

We say $\{x_n\}_{n=0}^k \subset (X,\kappa)$ is a {\em $\kappa$-path} (or a {\em path} if $\kappa$ is understood)
from $x_0$ to $x_k$ if $x_i \adjeq_{\kappa} x_{i+1}$ for $i \in \{0,\ldots,k-1\}$, and $k$ is the {\em length} of the path.

A subset $Y$ of a digital image $(X,\kappa)$ is
{\em $\kappa$-connected}~\cite{Rosenfeld},
or {\em connected} when $\kappa$
is understood, if for every pair of points $a,b \in Y$ there
exists a $\kappa$-path in $Y$ from $a$ to $b$.

\subsection{Digitally continuous functions}
Material in this section is largely quoted or paraphrased from~\cite{bs19a}.

We denote by $\id$ or $\id_X$ the
identity map $\id(x)=x$ for all $x \in X$.

\begin{definition}
{\rm \cite{Rosenfeld, Bx99}}
Let $(X,\kappa)$ and $(Y,\lambda)$ be digital
images. A function $f: X \to Y$ is 
{\em $(\kappa,\lambda)$-continuous}, or
{\em digitally continuous} when $\kappa$ and
$\lambda$ are understood, if for every
$\kappa$-connected subset $X'$ of $X$,
$f(X')$ is a $\lambda$-connected subset of $Y$.
If $(X,\kappa)=(Y,\lambda)$, we say a function
is {\em $\kappa$-continuous} to abbreviate
``$(\kappa,\kappa)$-continuous."
\end{definition}

\begin{thm}
{\rm \cite{Bx99}}
A function $f: X \to Y$ between digital images
$(X,\kappa)$ and $(Y,\lambda)$ is
$(\kappa,\lambda)$-continuous if and only if for
every $x,y \in X$, if $x \adj_{\kappa} y$ then
$f(x) \adjeq_{\lambda} f(y)$.
\end{thm}

\begin{thm}
\label{composition}
{\rm \cite{Bx99}}
Let $f: (X, \kappa) \to (Y, \lambda)$ and
$g: (Y, \lambda) \to (Z, \mu)$ be continuous 
functions between digital images. Then
$g \circ f: (X, \kappa) \to (Z, \mu)$ is continuous.
\end{thm}

\begin{definition}
Let $A \subset X$. A $\kappa$-continuous
function $r: X \to A$ is a {\em retraction}, and
{\em $A$ is a retract of $X$}, if $r(a)=a$ for
all $a \in A$. 
\end{definition}

A function $f: (X,\kappa) \to (Y,\lambda)$ is
an {\em isomorphism} (called a {\em homeomorphism}
in~\cite{Bx94}) if $f$ is a continuous bijection
such that $f^{-1}$ is continuous.

We use the following notation. For a
digital image $(X,\kappa)$,
\[ C(X,\kappa) = \{f: X \to X \, | \,
   f \mbox{ is $\kappa$-continuous}\}.
\]

Given $f \in C(X,\kappa)$, a point
$x \in X$ is a {\em fixed point of $f$} if
$f(x)=x$. We denote by $\Fix(f)$ the set
$\{x \in X \, | \, x 
   \mbox{ is a fixed point of } f \}$.
A point $x \in X$ is an {\em almost fixed point}~\cite{Rosenfeld,TsSm} or
an {\em approximate fixed point}~\cite{BEKLL} of $f$ if $x \adjeq_{\kappa} f(x)$.

We use the projection functions $p_1,p_2: \Z^2 \to \Z$
defined for $(x,y) \in \Z^2$ by $p_1(x,y) = x$,
$p_2(x,y)=y$. These functions are
$(c_1,c_1)$-continuous and $(c_2,c_1)$-continuous~\cite{Han05}.

\subsection{Freezing and cold sets}
Material in this section is largely quoted or paraphrased from~\cite{BxFpSets}.

Knowledge of $\Fix(f)$ for $f \in C(X,\kappa)$ can tell 
us much about $f|_{X \setminus \Fix(f)}$.
This motivates the study of freezing and cold sets.

\begin{definition}
\label{freezeDef}
{\rm \cite{BxFpSets}}
Let $(X,\kappa)$ be a digital image. We say $A \subset X$ is a 
{\em freezing set for $X$} if given $g \in C(X,\kappa)$,
$A \subset \Fix(g)$ implies $g=\id_X$. If no proper subset 
of a freezing set $A$ is a freezing set for $(X,\kappa)$, 
then $A$ is a {\em minimal freezing set}
\end{definition}

\begin{definition}
\label{s-cold-def}
{\rm \cite{BxFpSets}}
$A \subset X$ is a {\em cold set}
for the connected digital image $(X,\kappa)$
if given $g \in C(X,\kappa)$ such that
$g|_A = \id_A$, then for all $x \in X$, 
$g(x) \in N^*(X,x,\kappa)$
\end{definition}

\begin{remk}
{\rm \cite{BxFpSets}}
A freezing set is a cold set.
\end{remk}

\begin{defn}
\label{bdDef}
{\rm \cite{BxConvex}}
Let $X \subset \Z^n$.
\begin{itemize}
    \item The
{\em boundary of $X$ with respect to the $c_i$ adjacency},
$i \in \{1,2\}$, is
\[Bd_i(X) = \{x \in X \, | \mbox{ there exists } y \in \Z^n \setminus X \mbox{ such that } y \adj_{c_i} x\}.
\]
$Bd_1(X)$ is what is called the {\em boundary of $X$}
in~\cite{RosenfeldMAA}. This paper uses both $Bd_1(X)$ and $Bd_2(X)$.
\item The {\em interior of} $X$ with respect to the $c_i$ adjacency
is $Int_i(X) = X \setminus Bd_i(X)$.
\end{itemize}
\end{defn}

\begin{thm}
{\rm \cite{BxFpSets}}
\label{bdFreezes}
Let $X \subset \Z^n$ be finite. Then for $1 \le u \le n$, $Bd_1(X)$ is 
a freezing set for $(X,c_u)$.
\end{thm}

\begin{thm}
{\rm \cite{BxFpSets}}
\label{corners-min}
Let $X = \Pi_{i=1}^n [0,m_i]_{\Z}$.
Let $A = \Pi_{i=1}^n \{0,m_i\}$.
\begin{itemize}
\item Let $Y = \Pi_{i=1}^n [a_i,b_i]_{\Z}$ be
      such that $X \subset Y$. Let $f: X \to Y$ be
      $c_1$-continuous. If $A \subset \Fix(f)$, then $X \subset \Fix(f)$.
\item $A$ is a freezing set for $(X,c_1)$; minimal for $n \in \{1,2\}$.
\end{itemize}
\end{thm}

\begin{thm}
{\rm \cite{BxFpSets}}
\label{noProperSub}
Let $X = \prod_{i=1}^n[0,m_i]_{\Z}  \subset \Z^n$, where $m_i>1$ for all~$i$.
Then $Bd_1(X)$ is a minimal freezing set for $(X,c_n)$.
\end{thm}

\subsection{Digital disks and bounding curves}
Material in this section is largely quoted or paraphrased from~\cite{BxConvex}.

We say a finite $c_2$-connected 
set $S=\{x_i\}_{i=1}^n \subset \Z^2$ is a
{\em (digital) line segment} if the members of $S$ are collinear.

We say a segment with slope of $\pm 1$ is
{\em slanted}. An {\em axis-parallel} segment is horizontal or vertical.

\begin{remk}
\label{segSlope}
{\rm \cite{BxConvex}}
A digital line segment must be axis-parallel or slanted.
\end{remk}

A {\em closed curve} is a path $S=\{s_i\}_{i=0}^m$ 
such that $s_0=s_m$,
and $0 < |i - j| < m$ implies $s_i \neq s_j$. If
\[ N(S,x_0,\kappa) = N(S,x_m,\kappa) = \{x_1,x_{m-1}\} \mbox{ and }
\]
\[ 1 \le i < m \mbox{ implies } 
   N(S,x_i,\kappa) = \{x_{i-1}, x_{i+1}\},
\]
$S$ is a {\em cycle}. We may also refer to a cycle as
a {\em (digital) $\kappa$-simple closed curve}.
For a simple closed curve $S \subset \Z^2$ we generally assume
\begin{itemize}
    \item $m \ge 8$ if $\kappa = c_1$, and
    \item $m \ge 4$ if $\kappa = c_2$.
\end{itemize}
These requirements are necessary for the Jordan Curve
Theorem of digital topology, below, as a
$c_1$-simple closed curve in $\Z^2$ must have at least 8 points to
have a nonempty finite complementary $c_2$-component,
and a $c_2$-simple closed curve in $\Z^2$ must have at least 4 points to
have a nonempty finite complementary $c_1$-component.
Examples in~\cite{RosenfeldMAA} show why it is
desirable to consider $S$ and $\Z^2 \setminus S$
with different adjacencies.

\begin{thm}
{\rm \cite{RosenfeldMAA}}
{\em (Jordan Curve Theorem for digital topology)}
Let $\{\kappa, \kappa'\} = \{c_1, c_2\}$.
Let $S \subset \Z^2$ be a simple closed 
$\kappa$-curve such that $S$ has at least 8 points if
$\kappa = c_1$ and such that $S$ has at least 
4 points if $\kappa = c_2$. Then
$\Z^2 \setminus S$ has exactly 2 $\kappa'$-connected
components.
\end{thm}

One of the $\kappa'$-components of 
$\Z^2 \setminus S$ is finite and the other is infinite. This 
suggests the following.
\begin{defn}
\label{diskDef}
{\rm \cite{BxConvex}}
Let $S \subset \Z^2$ be a $c_2$-closed curve such that
$\Z^2 \setminus S$ has two $c_1$-components, one finite and the
other infinite. The union $D$ of $S$ and the finite $c_1$-component 
of $\Z^2 \setminus S$ is a {\em (digital) disk}. $S$ is
a {\em bounding curve} of $D$. The finite $c_1$-component 
of $\Z^2 \setminus S$ is the {\em interior of} $S$, denoted $Int(S)$,
and the infinite $c_1$-component of $\Z^2 \setminus S$ is the {\em exterior of} 
$S$, denoted $Ext(S)$.
\end{defn}

\begin{figure}
    \centering
    \includegraphics[height=1.5in]{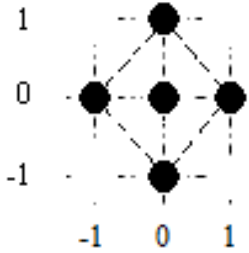}
    \caption{\cite{BxConvex} The $c_1$-disk
    $D = \{(x,y) \in \Z^2 \, | \, |x| + |y| < 2\}$.
    The bounding curve $S = \{(x,y) \in \Z^2 \, | \,
    |x| + |y| =1\} = D \setminus \{(0,0)\}$ is not $c_1$-connected. \newline
    }
    \label{fig:diamond}
\end{figure}

\begin{figure}
    \centering
    \includegraphics[height=1.25in]{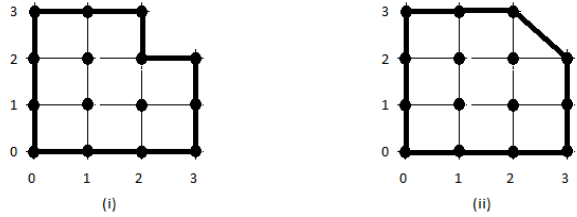}
    \caption{\cite{BxConvex} Two views of $D = [0,3]_{\Z}^2 \setminus \{(3,3)\}$, 
    which can be regarded as a $c_1$-disk with either of the
    closed curves shown in dark as a bounding curve. \newline
    (i) The dark line segments show a $c_1$-simple closed curve $S$
    that is a bounding curve for~$D$.
    Note the point $(2,2)$ in the bounding curve shown.
    By Definition~\ref{bdDef},
    $(2,2) \not \in Bd_1(D)$; however, $(2,2) \in Bd_2(D)$. \newline
    (ii) The dark line segments show a $c_2$-closed curve $S$
    that is a minimal bounding curve for~$D$.
       }
    \label{fig:2sccBdry}
\end{figure}

\begin{figure}
    \centering
    \includegraphics{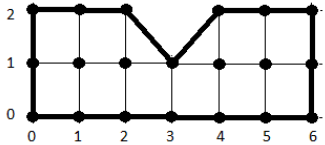}
    \caption{\cite{BxConvex} $D = [0,6]_{\Z} \times [0,2]_{\Z} \setminus \{(3,2)\}$
     shown with a bounding curve $S$ in dark segments. $D$ is
     not a disk with either the $c_1$ or the $c_2$ adjacency,
     since with either of these adjacencies,
     $\Z^2 \setminus S$ has two bounded components,
     $\{(1,1), (2,1)\}$ and $\{(4,1), (5,1)\}$.
     }
    \label{fig:notDisk}
\end{figure}

Notes:
\begin{itemize}
    \item If $D$ is a digital disk determined as above by a bounding 
    $c_2$-closed curve $S$, then $(S,c_1)$ can be 
    disconnected. See Figure~\ref{fig:diamond}.
    \item There may be more than one closed curve $S$
          bounding a given disk~$D$. See Figure~\ref{fig:2sccBdry}.
          When $S$ is understood as a bounding curve of a disk $D$,
          we use the notations $Int(S)$ and $Int(D)$ interchangeably.
    \item Since we are interested in finding {\em minimal}
          freezing or cold sets and since it turns out we often compute these
          from bounding curves, we may prefer those of
          minimal size. A bounding curve~$S$
          for a disk $D$ is {\em minimal} if there is no
          bounding curve $S'$ for $D$ such that
          $\#S' < \#S$.
    \item In particular, a bounding 
          curve need not be contained in $Bd_1(D)$.
          E.g., in the disk~$D$
          shown in Figure~\ref{fig:2sccBdry}(i), $(2,2)$ is a point
          of the bounding curve; however, all of the points
          $c_1$-adjacent to $(2,2)$ are members of~$D$, so
          by Definition~\ref{bdDef}, $(2,2) \not \in Bd_1(D)$.
          However, a bounding curve for $D$ must be contained
          in $Bd_2(D)$.
    \item In Definition~\ref{diskDef}, we use $c_2$ adjacency for
          $S$ and we do not require $S$ to be simple. 
          Figure~\ref{fig:2sccBdry} shows why these seem
          appropriate.
          \begin{itemize}
              \item The $c_2$ adjacency allows slanted
              segments in bounding curves and makes possible
          a bounding curve in subfigure~(ii) with fewer points
          than the bounding curve in subfigure~(i) in which
          adjacent pairs of the bounding curve are restricted
          to $c_1$ adjacency.
          \item Neither of the bounding curves shown in
                Figure~\ref{fig:2sccBdry} is a $c_2$-simple closed
                curve. E.g., non-consecutive points of each of
                the bounding curves,
                $(0,1)$ and $(1,0)$, are $c_2$-adjacent. The
                bounding curve shown in 
                Figure~\ref{fig:2sccBdry}(ii) is clearly also not a
                $c_1$-simple closed curve.
          \end{itemize}
    \item A closed curve that is not simple may be the boundary~$Bd_2$
          of a digital image that is not a disk. This is illustrated
          in Figure~\ref{fig:notDisk}.
\end{itemize}

More generally, we have the following.
\begin{defn}
\label{boundingCurvesDef}
{\rm \cite{BxConvex}}
Let $X \subset \Z^2$ be a finite, $c_i$-connected set,
$i \in \{1,2\}$. Suppose there are 
pairwise disjoint $c_2$-closed curves
$S_j \subset X$, $1 \le j \le n$, such that
\begin{itemize}
    \item $X \subset S_1 \cup Int(S_1)$;
    \item for $j>1$, $D_j = S_j \cup Int(S_j)$ is a digital disk;
    \item no two of 
    \[ S_1 \cup Ext(S_1), D_2, \ldots, D_n
    \]
    are $c_1$-adjacent or $c_2$-adjacent; and
    \item we have 
    \[ \Z^2 \setminus X = Ext(S_1) \cup \bigcup_{j=2}^n Int(S_j).
    \]
\end{itemize}
Then $\{S_j\}_{j=1}^n$ is a {\em set of bounding curves of} $X$.
\end{defn}

Note: As above, a digital image $X \subset \Z^2$ may have more than one
set of bounding curves.

 \begin{figure}
        \centering
        \includegraphics[height=1in]{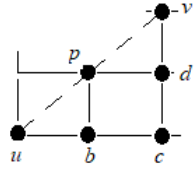}
        \caption{\cite{BxConvex}
        $p \in \overline{uv}$ in a bounding curve,
        with $\overline{uv}$ slanted.
        Note $u \not \adj_{c_1} p \not \adj_{c_1} v$,
        $p \adj_{c_2} c \not \adj_{c_1} p$,
        $\{p,c\} \subset N(\Z^2,c_1,b) \cap N(\Z^2,c_1,d)$. If
        $X$ is slant-thick at $p$ then $c \in X$.
        (Not meant to be understood as showing all of $X$.)}
        \label{fig:innerBdPt}
    \end{figure}

\subsection{Thickness}
A notion of ``thickness" in a digital image $X$, introduced in~\cite{BxConvex},
means, roughly speaking, $X$ is ``locally" like a disk.

Our definition of thickness depends on a notion of an
``interior angle" of a disk. We have the following.

\begin{definition}
\label{interiorAngleDef}
{\rm \cite{BxConvex}}
Let $s_1$ and $s_2$ be sides of a digital disk
$X \subset \Z^2$, i.e., maximal digital line segments
in a bounding curve $S$ of $X$, such that 
$s_1 \cap s_2 = \{p\} \subset X$.
The {\em interior angle of $X$ at $p$} is the
angle formed by $s_1$, $s_2$, and $Int(S)$.
\end{definition}

\begin{figure}
        \centering
        \includegraphics[height=1in]{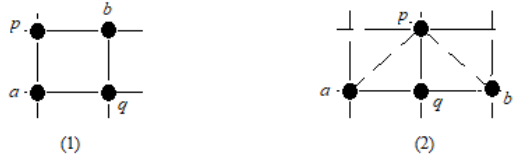}
        \caption{\cite{BxConvex} 
        (1) $\angle apb$ is a
        $90^{\circ}$ ($\pi/2$ radians)
        angle of a bounding curve of $X$ at $p \in A_1$, with
        horizontal and vertical sides. If $X$ is $90^{\circ}$-thick at $p$ then
        $q \in Int(X)$. (Not meant to
        be understood as showing all of $X$.)
        \newline
        (2) $\angle apb$ is a
        $90^\circ$ ($\pi/2$ radians) angle
         between slanted segments of a bounding curve. If $X$ is
         $90^{\circ}$-thick at $p$ then $q \in Int(X)$ and therefore
         $q' \in X$. (Not meant to
        be understood as showing all of $X$).
}
        \label{fig:degrees90a}
    \end{figure}

\begin{figure}
        \centering
        \includegraphics[height=1in]{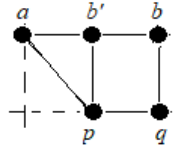}
        \caption{\cite{BxConvex} $\angle apq$ is an angle of
        135$^ \circ$ degrees ($3 \pi /4$ radians)
        of a bounding curve of $X$ at $p$, with
        $\overline{ap} \cup \overline{pq}$
            a subset of the bounding curve. If
            $X$ is $135^{\circ}$-thick at $p$ then $b,b' \in X$. (Not meant to
        be understood as showing all of $X$.)
        }
        \label{fig:degrees135c1}
        \end{figure}

\begin{definition}
{\rm \cite{BxConvex}}
Let $X \subset \Z^2$ be a digital disk. Let $S$ be a bounding curve of $X$ and $p \in S$.
\begin{itemize}
    \item Suppose $p$ is in a maximal slanted segment $\sigma$ of $S$ such that
          $p$ is not an endpoint of $\sigma$. Then $X$ is {\em slant-thick at $p$} 
          if there exists $c \in X$ such that (see Figure~\ref{fig:innerBdPt})
          \begin{equation}
          \label{slantSegProp}
           c \adj_{c_2} p \not \adj_{c_1} c,
          \end{equation}
    \item Suppose $p$ is the vertex of a 90$^\circ$ ($ \pi / 2$ radians) interior angle
              $\theta$ of $S$. Then $X$ is {\em $90^\circ$-thick at $p$} if there
              exists $q \in Int(X)$ such that
              \begin{itemize}
                  \item if $\theta$ has axis-parallel sides then
                        $q \adj_{c_2} p \not \adj_{c_1} q$ (see Figure~\ref{fig:degrees90a}(1));
                  \item if $\theta$ has slanted sides then $q \adj_{c_1} p$ (see
                        Figure~\ref{fig:degrees90a}(2)).
              \end{itemize}
    \item Suppose $p$ is the vertex of a 135$^\circ$ ($3 \pi / 4$
          radians) interior angle $\theta$ of $S$. Then $X$ is
          {\em 135$^\circ$-thick at $p$} if there exist $b,b' \in X$
          such that $b$ and $b'$ are in the interior of $\theta$ and
          (see Figure~\ref{fig:degrees135c1})
          \[ b \adj_{c_2} p \not \adj_{c_1} b~~~ \mbox{ and }~~~ 
                b' \adj_{c_1} p.
          \]
\end{itemize}
\end{definition}

\begin{defn}
\label{thickness}
{\rm \cite{BxConvex,BxColdAndFreezing}}
Let $X \subset \Z^2$ be a digital disk. We say $X$ is
{\em thick} if the following are satisfied. For some bounding
curve $S$ of $X$,
\begin{itemize}
    \item for every maximal slanted segment of~$S$,
          if  $p \in S$ is not an endpoint of  $S$, 
          then $X$ is slant-thick at $p$, 
          and
    \item for every $p$ that is the vertex of a 90$^\circ$ ($ \pi / 2$
          radians) interior angle $\theta$ of $S$, $X$ is {\em $90^\circ$-thick at $p$}, 
          and
\item for every $p$ that is the vertex of a 135$^\circ$ ($3 \pi / 4$
      radians) interior angle $\theta$ of $S$, $X$ is 135$^\circ$-thick at $p$.
\end{itemize}
\end{defn}

\subsection{Convexity}
A set $X$ in a Euclidean space $\R^n$ is
{\em convex} if for every pair of distinct
points $x,y \in X$, the line segment
$\overline{xy}$ from $x$ to $y$ is contained in $X$.
The {\em convex hull of} $Y \subset \R^n$,
denoted $hull(Y)$, is the
smallest convex subset of $\R^n$ that contains~$Y$.
If $Y \subset \R^2$ is a finite set, then
$hull(Y)$ is a single point if $Y$ is a singleton;
a line segment if $Y$ has at least 2 members and all are
collinear; otherwise, $hull(Y)$ is a polygonal disk,
and the endpoints of the edges of $hull(Y)$ are its {\em vertices}.

A digital version of convexity can be stated
for subsets of the digital plane~$\Z^2$ as follows.
A finite set $Y \subset \Z^2$ is 
{\em (digitally) convex} {\rm \cite{BxConvex}} if either
\begin{itemize}
    \item $Y$ is a single point, or
    \item $Y$ is a digital line segment, or
    \item $Y$ is a digital disk with a bounding curve $S$
          such that the endpoints of the maximal line segments
          of~$S$ are the vertices of $hull(Y) \subset \R^2$.
\end{itemize}

\section{Tools for determining fixed point sets}
The following assertions will be useful in determining fixed point and freezing
sets.

\begin{prop}
{\rm (Corollary 8.4 of~\cite{bs19a})}
\label{uniqueShortestProp}
Let $(X,\kappa)$ be a digital image and
$f \in C(X,\kappa)$. Suppose
$x,x' \in \Fix(f)$ are such that
there is a unique shortest
$\kappa$-path $P$ in~$X$ from $x$ 
to $x'$. Then $P \subseteq \Fix(f)$.
\end{prop}

Lemma~\ref{c1pulling}, below,
\begin{quote}
... can
be interpreted to say that
in a $c_u$-adjacency,
a continuous function that
moves a point~$p$ also moves
a point that is ``behind"
$p$. E.g., in $\Z^2$, if $q$ and $q'$ are
$c_1$- or $c_2$-adjacent with $q$
left, right, above, or below $q'$, and a
continuous function $f$ moves $q$ to the left,
right, higher, or lower, respectively, then
$f$ also moves $q'$ to the left,
right, higher, or lower, respectively~\cite{BxFpSets}.
\end{quote}

\begin{lem}
\label{c1pulling}
{\rm ~\cite{BxFpSets}}
Let $(X,c_u)\subset \Z^n$ be a digital image, 
$1 \le u \le n$. Let $q, q' \in X$ be such that
$q \adj_{c_u} q'$.
Let $f \in C(X,c_u)$.
\begin{enumerate}
    \item If $p_i(f(q)) > p_i(q) > p_i(q')$
          then $p_i(f(q')) > p_i(q')$.
    \item If $p_i(f(q)) < p_i(q) < p_i(q')$
          then $p_i(f(q')) < p_i(q')$.
\end{enumerate}
\end{lem}

\begin{remk}
\label{boundingCurveSetFreezes}
{\rm \cite{BxFpSets}}
If $X \subset \Z^2$ is finite, then
a set of bounding curves for $X$ is a freezing set
for $(X,c_i)$, $i \in \{1,2\}$.
\end{remk}

In particular, we have:

\begin{thm}
\label{bdCurveFreezes}
Let $D$ be a digital disk in $\Z^2$. Let
$S$ be a bounding curve for $D$. Then $S$ is
a freezing set for $(D,c_1)$ and for $(D,c_2)$.
\end{thm}

The next two results form a dual pair.

\begin{thm}
\label{convDiskThmActual}
{\rm \cite{BxConvex}}
Let $X$ be a thick convex disk with a
    bounding curve $S$.
    Let $A_1$ be the set of points $x \in S$ such that
$x$ is an endpoint of a maximal axis-parallel edge of $S$. Let $A_2$ 
be the union of slanted line segments in $S$.
Then $A = A_1 \cup A_2$ is a minimal 
freezing set for $(X,c_1)$.
\end{thm}

\begin{thm}
\label{convDiskThmC2Actual}
{\rm \cite{BxConvex}}
Let $X$ be a thick convex disk with a minimal bounding
curve $S$. Let $B_1$ be the set of
points $x \in S$ such that
$x$ is an endpoint of a maximal
slanted edge in $S$. Let $B_2$ 
be the union of maximal axis-parallel line segments in $S$.
Let $B = B_1 \cup B_2$. Then $B$ is a 
minimal freezing set for $(X,c_2)$.
\end{thm}

The next two results form another dual pair, generalizing the
previous pair.

\begin{thm}
{\rm \cite{BxSubsets}}
\label{corners}
Let $V_i \subset X \subset \Z^2$, $i \in \{1,\ldots,n\}$ 
where each $V_i$ is a thick convex disk. 
Let $X' = \bigcup_{i=1}^n V_i$.
Let $C_i$ be a bounding curve of $V_i$.
Let $A_{1,i}$ be the set of endpoints of maximal
horizontal or vertical segments of $C_i$.
Let $A_{2,i}$ be the union of maximal slanted
segments of $C_i$. Then
$A = (X \setminus X') \cup \bigcup_{i=1}^n (A_{1,i} \cup A_{2,i})$ is a freezing set for $(X,c_1)$.
\end{thm}

\begin{thm}
{\rm \cite{BxSubsets}}
\label{slantCorners}
Let $V_i \subset X \subset \Z^2$, $i \in \{1,\ldots,n\}$ 
where each $V_i$ is a thick convex disk. 
Let $X' = \bigcup_{i=1}^n V_i$.
Let $C_i$ be a bounding curve of $V_i$.
Let $B_{1,i}$ be the union of maximal horizontal 
and maximal vertical segments of $C_i$.
Let $B_{2,i}$ be the set of endpoints of maximal slanted
segments of $C_i$. Then
$B = (X \setminus X') \cup \bigcup_{i=1}^n (B_{1,i} \cup B_{2,i})$ is a freezing set for $(X,c_2)$ (the adjacency
is misprinted as $c_1$ in~{\rm \cite{BxSubsets}}).
\end{thm}

\section{Unifying sets}
\subsection{Definition and general properties}
\begin{definition}
\label{harmonyDef}
Let $(X,\kappa)$ be a digital image. Let $A \subset X$.
Suppose whenever $f,g \in C(X,\kappa)$ are such that
$f(A)=g(A)=A$ and
$f|_A = g|_A$, we have $f=g$. Then we say $A$ is a
{\em unifying set} for $(X,\kappa)$. $A$ is a
{\em minimal unifying set} if $A$ is a unifying set
and no proper subset of $A$ is a unifying set for
$(X,\kappa)$.
\end{definition}

\begin{remk}
\label{unifyImpliesFreeze}
Observe:
\begin{itemize}
    \item By taking $g$ to be the identity function 
          $\id_X$ in Definition~\ref{harmonyDef}, 
          we see that a unifying set is a freezing set.
          We have not determined whether the converse
          is true.
    \item It is trivial that $X$ is a unifying set for
          $(X,\kappa)$. We are therefore interested
          in finding minimal unifying sets. In light
          of the above, a minimal freezing set is a ``good
          candidate" for a minimal unifying set.
\end{itemize}
\end{remk}

In the following, we study conditions
for which a freezing set must be unifying.

The desirability of the requirement that $f(A)=g(A)=A$ in 
Definition~\ref{harmonyDef} is illustrated in the 
following, in which this requirement is not met.

\begin{exl}
Let $X = [0,m]_{\Z} \times [0,n]_{\Z}$ for $m \ge 2,~n >0$.
Let $f,g: X \to X$ be the functions
\[ f(x,y) = (0,y), ~~~~~
g(x,y) = \left \{ \begin{array}{ll}
             (0,y) & \mbox{if } x \in \{0,m\}; \\
             (1,y) & \mbox{if } 1 \le x \le m-1,
            \end{array} \right .
\]
We take
\[ A = \{(0,0), (0,n), (m,0), (m,n)\}.
\]
Note by Theorem~\ref{corners-min}, $A$ is a minimal freezing 
set for $(X,c_1)$. We see easily that $f,g \in C(X,c_1)$,
$f|_A = g|_A$, $f(A)=g(A)$ is a proper subset of $A$,
and $f \neq g$.
\end{exl}

The following shows that unifying sets are preserved by
isomorphism.

\begin{thm}
\label{isoInvariant}
Let $(X,\kappa)$ and $(Y,\lambda)$ be digital images such that
there exists an isomorphism $F: (X,\kappa) \to (Y,\lambda)$.
If $A$ is a unifying set for $(X,\kappa)$ then
$F(A)$ is a unifying set for $(Y, \lambda)$.
\end{thm}

\begin{proof}
Let $f,g \in C(Y,\lambda)$ such that
$f(F(A))=g(F(A)) = F(A)$ and $f|_{F(A)} = g|_{F(A)}$.

We have, by Theorem~\ref{composition}, 
$f' = F^{-1} \circ f \circ F, ~ g'= F^{-1} \circ g \circ F \in C(X,\kappa)$,
and for $a \in A$ we have $f \circ F(a) = g \circ F(a)$, so
\[ f'(a)= F^{-1} \circ f \circ F(a) = F^{-1} \circ g \circ F(a),
   =g'(a).
\]
Also, given $b = F(a)$ for $a \in A$, by assumption we have
$f(b)=g(b)$, hence 
\[ f'(a) = F^{-1}(f(b)) = F^{-1}(g(b)) = g'(a).
\]
Since $A$ is unifying, $f' = g'$. Therefore,
\[ f = F \circ f' \circ F^{-1} = F \circ g' \circ F^{-1} = g,
\]
so $F(A)$ is unifying for $(Y,\lambda)$.
\end{proof}

We have the following generalization of
Proposition~\ref{uniqueShortestProp}.

\begin{prop}
\label{uniqueCommonShortest}
Let $f,g: X \to Y$ such that $f$ and $g$ are both
$(\kappa,\lambda)$-continuous. Suppose $x_0,x_1 \in X$
and there is a $\kappa$-path $P$ of length $n$ in $X$ from
$x_0$ to $x_1$. Suppose $y_0=f(x_0)=g(x_0)$,
$y_1=f(x_1)=g(x_1)$, and there is a unique shortest
path $Q$ of length $n$ in $Y$ from $y_0$ to $y_1$. Then 
$f(P)=g(P)=Q$ and $f|_P = g|_P$.
\end{prop}

\begin{proof}
Since $f(P)$ and $g(P)$ must be $\lambda$-paths
from $y_0$ to $y_1$, our uniqueness and length 
restrictions imply $f(P)=g(P)=Q$. Continuity
implies $f|_P = g|_P$.
\end{proof}

\subsection{Cycles}
\label{unifyingCycle}

\begin{thm}
\label{3ptsForCycles}
{\rm \cite{BxFpSets}}
Let $n > 4$. Let $C_n = \{x_m\}_{m=0}^{n-1} \subset \Z^2$, 
where the members of $C_n$ are indexed circularly.
Let $A = \{x_i,x_j,x_k\}$ be a set of
distinct members of $C_n$ such
that $C_n$ is a union of unique shorter paths
determined by these points.
Then $A$ is a minimal freezing set for $C_n$.
\end{thm}

\begin{thm}
The set $A$ of Theorem~\ref{3ptsForCycles} is
a unifying set for $(C_n,\kappa)$, and any
$f \in C(X,\kappa)$ such that $f(A)=A$ must be an
isomorphism of $(X,\kappa)$.
\end{thm}

\begin{proof}
Let $\widehat{x_ix_j}$, $\widehat{x_ix_k}$, and
$\widehat{x_jx_k}$ be the unique shorter paths in
$C_n$ from $x_i$ to $x_j$, from $x_i$ to $x_k$, and
from $x_j$ to $x_k$, respectively. Let
$B = \{ \widehat{x_ix_j}, \widehat{x_ix_k}, \widehat{x_jx_k}\}$. 
Let $f,g \in C(C_n,\kappa)$ such that 
\begin{equation} 
\label{arcPreserved}
f(A)=g(A)=A \mbox{ and } f|_A = g|_A.
\end{equation}

Suppose $f \neq g$. Consider the following cases.
\begin{itemize}
    \item The members of $B$ have distinct lengths.
          Without loss of generality,
          \begin{equation}
          \label{orderLengths}
          length(\widehat{x_ix_j}) <
           length(\widehat{x_ix_k}) < 
           length(\widehat{x_jx_k}).
           \end{equation}
          Since we have that
          both $f(\widehat{x_ix_j})$ and $g(\widehat{x_ix_j})$ are paths of length
          at most $length(\widehat{x_ix_j})$ from
          $f(x_i)=g(x_i)$ to $f(x_j)=g(x_j)$, 
          from~(\ref{orderLengths}) and
          Proposition~\ref{uniqueCommonShortest},
          $f(\widehat{x_ix_j})=g(\widehat{x_ix_j})$ and
          $f|_{\widehat{x_ix_j}} = g|_{\widehat{x_ix_j}}$
          is a bijection of $\widehat{x_ix_j}$.
          Indeed, we must have that $f$ and $g$
          coincide with $\id_X$ on $\widehat{x_ix_j}$,
          for otherwise we would have
          $f(x_i)=g(x_i)=x_j$, $f(x_j)=g(x_j)=x_i$,
          $f(x_k)=g(x_k)=x_k$, so $f(\widehat{x_ix_k})$
          is a $\kappa$-path from $x_j$ to $x_k$,
          contrary to~(\ref{orderLengths}). Then
          by~(\ref{arcPreserved}) we have 
          $f|_A = g|_A = \id_A$, and from
          Proposition~\ref{uniqueCommonShortest}
          it follows that $f=g=\id_X$.
    \item Suppose two members, but not all three,
          of $B$ have the same
          length; without loss of generality,
          $length(\widehat{x_ix_j}) =
           length(\widehat{x_ix_k})$. Then
           either $f|_A = g|_A = \id_A$ or
           $f(x_i)=g(x_i)=x_i$, $f(x_j)=g(x_j)=x_k$, and
           $f(x_k)=g(x_k)=x_j$. Then much as above,
           $f = g$ is an isomorphism of $(X,\kappa)$.
    \item Suppose all three members of $B$ have the same
          length. Then $f|_A = g|_A$ is a permutation
          of $A$. Much as above, it follows that
          $f = g$ is an isomorphism of $(X,\kappa)$.
\end{itemize}
In all cases, we concluded that $f = g$ is an 
isomorphism of $(X,\kappa)$. Thus $A$ is a unifying set
for $(X,\kappa)$.
\end{proof}

\subsection{Trees}
A {\em tree} is a connected acyclic graph $(X,\kappa)$.
By {\em acyclic} we mean lacking any 
closed curve of more than 2 points.
The {\em degree} of a vertex $x$ in $X$ is the number 
of distinct vertices $y \in X$ such that $x \adj y$.

\begin{thm}
{\rm \cite{BxFpSets}}
\label{tree}
Let $(X,\kappa)$ be a digital image such that
the graph $G=(X,\kappa)$ is a finite tree with $\#X > 1$.
Let $A$ be the set of vertices of $G$ that have degree 1. 
Then $A$ is a minimal freezing set for $G$.
\end{thm}

\begin{thm}
Let $(X,\kappa)$ be a digital image such that
the graph $G=(X,\kappa)$ is a finite tree with $\#X > 1$.
Let $A$ be the set of vertices of $G$ that have degree~1. 
Then $A$ is a minimal unifying set for $G$. Also, if
$f \in C(X,\kappa)$ such that $f(A)=A$, then
$f$ is an isomorphism of $(X,\kappa)$.
\end{thm}

\begin{proof}
Let $a_0 \in A$. Since $X$ is finite, we have that 
$A$ is also finite - say, $A=\{a_i\}_{i=0}^n$. Since
$G$ is a tree, for $0< i \le n$ there is a unique
shortest $\kappa$-path $P_i$ in $X$ from $a_0$ to $a_i$.
Let $L = \{\ell_j\}_{j=1}^m$ be the set of
distinct lengths of the members of $\{P_i\}_{i=1}^n$, with
\[ \ell_1 < \ell_2 < \ldots < \ell_m.
\]
Let $L_j = \{P_i ~|~ length(P_i) = \ell_j\}.$
Since $A$ is finite and $f(A)=g(A)=A$, 
\begin{equation}
    \label{bijectA}
    f|_A = g|_A: A \to A \mbox{ is a bijection.}
\end{equation}

Let $f,g \in C(X,\kappa)$ be such that
$f(A)=g(A)=A$ and $f|_A = g|_A$. Every $P_k$ of
length $\ell_1$ is the unique shortest $\kappa$-path in X
from $a_0$ to some $a_k \in A \setminus \{a_0\}$.
Since $f(P_k)$ is a path from $f(a_0)=g(a_0)$ to
$f(a_k)=g(a_k)$, our choice of $\ell_1$ 
and Proposition~\ref{uniqueCommonShortest} imply
$f|_{P_k} = g|_{P_k}$, $f(P_k) = g(P_k)$ has
length $\ell_1$, and from~(\ref{bijectA}) that
$f|_{L_1} = g|_{L_1}$ is a bijection of $L_1$.
It follows easily that $f|_{L_1} = g|_{L_1}$ is an
isomorphism. This provides the base case
of an induction argument.

Suppose $u \in \Z$, $0 \le u < m$;
$f|_{P_k} = g|_{P_k}$ for every 
$P_k \in \bigcup_{j=1}^u L_j$; and
\begin{equation}
\label{bijectArcs}
   f|_{\bigcup_{j=1}^u L_j}=g|_{\bigcup_{j=1}^u L_j}
\mbox{ is a bijection of } \bigcup_{j=1}^u L_j.
\end{equation}
Now consider
$P_k \in L_{u+1}$. $f(P_k)$ and $g(P_k)$ are $\kappa$-paths
in $X$ from $f(a_0)=g(a_0)$ to $f(a_k)=g(a_k)$
of length at most $\ell_{u+1}$. 
By~(\ref{bijectArcs}), $f(P_k)$ and $g(P_k)$ cannot have
length less than $\ell_{u+1}$. Therefore,
each of $f(P_k)$ and $g(P_k)$ belongs to $L_{u+1}$.
By the uniqueness condition that defines $L_{u+1}$
it follows that $f|_{P_k}=g|_{P_k}$.
By~(\ref{bijectA}), $f|_{L_{u+1}}=g|_{L_{u+1}}$ is
a bijection. It follows from the above that
$f|_{\bigcup_{j=1}^{u+1} L_j}=g|_{\bigcup_{j=1}^{u+1} L_j}$
is a bijection of $\bigcup_{j=1}^{u+1} L_j$, and, further,
an isomorphism.

This completes the induction. Since 
$X=\bigcup_{j=1}^m L_j$, we have
$f = g$. Since $f$ was chosen arbitrarily,
$A$ is a unifying set. Also, $f$ is an isomorphism.

To show the minimality of $A$, we see easily that
for any $a \in A$ there is a $\kappa$-retraction
$r: X \to X \setminus \{a\}$, so $r$ and $\id_X$
are members of $C(X,\kappa)$ that coincide on
$A \setminus \{a\}$, 
$r(A \setminus \{a\}) = \id_X(A \setminus \{a\}) =
(A \setminus \{a\})$, but $r \neq \id_X$.
\end{proof}

\subsection{Complete graphs}
\begin{thm}
Let $(X,\kappa)$ be a digital image that is a 
complete graph, where $\#X > 1$. Let $A \subset X$.
Then the following are equivalent.
\begin{enumerate}
    \item $A = X$.
    \item $A$ is a unifying set for $(X,\kappa)$.
    \item $A$ is a freezing set for $(X,\kappa)$.
\end{enumerate}
\end{thm}

\begin{proof}
$1) \Rightarrow 2) \Rightarrow 3)$: These implications
are noted in Remark~\ref{unifyImpliesFreeze}.

$3) \Rightarrow 1)$: Suppose otherwise. Then there exists
$x_0 \in X \setminus A$. Let $x_1 \in X \setminus \{x_0\}$.
Let $g: X \to X$ be defined by
\[ g(x) = \left \{ \begin{array}{ll}
     x & \mbox{for } x \neq x_0; \\
     x_1 & \mbox{for } x = x_0.
     \end{array} \right .
\]
Since $(X,\kappa)$ is a complete graph, 
$g \in C(X,\kappa)$. Note 
$g|_A = \id_A$. But since $g(x_0) \neq x_0$, 
we have a contradiction of the assumption that $A$ 
is freezing. The contradiction gives us the 
desired conclusion.
\end{proof}

\subsection{Rectangles in $\Z^2$ with axis-parallel sides and $c_1$}
\label{c1RectangleSec}
In this section, we study unifying sets for
digital rectangles with axis-parallel edges in~$\Z^2$,
using the $c_1$ adjacency.

\begin{prop}
\label{essentialc2cold90}
{\rm \cite{BxColdAndFreezing}}
Let $X \subset \Z^2$. Let $S$ be a minimal bounding curve for $X$. 
Let $p_0$ be the vertex of an interior angle of $S$, formed 
by axis-parallel edges $E_1$ and $E_2$ of $S$, 
of measure $90^{\circ}$ ($\pi / 2$ radians). Let $A$ be any of a freezing set for $(X,c_1)$, 
a cold set for $(X,c_1)$, a freezing set for $(X,c_2)$, or a cold set for $(X,c_2)$.
Let $X$ be $90^{\circ}$-thick at $p_0$. Then $p_0 \in A$.
\end{prop}

\begin{prop}
\label{c1minRectFreeze}
Let $m>1$, $n>1$, and $X = [0,m]_{\Z} \times [0,n]_{\Z}$. Let
$A \subset X$. Then $A$ is a freezing set for
$(X,c_1)$ if and only if
\[ A' = \{(0,0), (m,0), (0,n), (m,n)\} \subset A.
\]
Therefore, $A'$ is the only minimal freezing set
for $(X,c_1)$.
\end{prop}

\begin{proof}
If $A$ is a freezing set, then by 
Proposition~\ref{essentialc2cold90}, $A' \subset A$.
Since $A'$ is a freezing set by Theorem~\ref{corners},
it follows that $A'$ is unique as a minimal freezing set.

If $A' \subset A$ then, since $A'$ is a freezing set,
$A$ is a freezing set~\cite{BxFpSets}.
\end{proof}

\begin{thm}
\label{axisParallelUnifying}
Let $X = [-m,m]_{\Z} \times [-n,n]_{\Z}$. Let 
\[A = \{(-m,-n), (-m,n), (m, -n), (m,n)\}.
\]
Then $A$ is a unifying set for $(X,c_1)$. Further,
every $f \in C(X,c_1)$ such that $f(A)=A$ is an
isomorphism.
\end{thm}

\begin{proof}
Let $f,g \in C(X,c_1)$ be such that
$f(A)=g(A)=A$ and $f|_A = g|_A$. 
Let $B,T,L,R$ be the bottom, top, left, and right
edges, respectively:
\[ B = [-m,m]_{\Z} \times \{-n\},~~~~~
   T = [-m,m]_{\Z} \times \{n\},
\]
\[ L = \{-m\} \times [-n,n]_{\Z},~~~~~
   R = \{m\} \times [-n,n]_{\Z}.
\]
Consider the following cases.
\begin{itemize}
    \item $m<n$. Since $f(A)=g(A)=A$, we have
          that $f(B)$, $g(B)$, $f(T)$, and $g(T)$
          are $c_1$-paths of length at most 2m between 
          distinct members of $A$, and since the closest
          distinct members of $A$ are joined by paths
          of length $2m$, $f(B)$, $g(B)$, $f(T)$,
          and $g(T)$ are paths of length $2m$. Therefore, 
          $f(B \cup T)= g(B \cup T)=B \cup T$.
          Continuity implies that for all
          $(x,y)\in B \cup T$,
          one of the following holds:
          \begin{itemize}
              \item $f(x,y)=g(x,y)=(x,y)$, or
              \item $f(x,y)=g(x,y)=(-x,y)$, or
              \item $f(x,y)=g(x,y)=(x,-y)$, or
              \item $f(x,y)=g(x,y)=(-x,-y)$.
          \end{itemize}
          
          Suppose the first case, 
          $f(x,y)=g(x,y)=(x,y)$ for $(x,y)\in B \cup T$.
          Each $(x,y) \in X$ lies on the unique shortest
          $c_1$-path between $b=(x,-n)$ and $t=(x,n)$.
          Since $f(b)=g(b)=b$ and $f(t)=g(t)=t$, we must
          have $f(x,y)=g(x,y)=(x,y)$ by
          Proposition~\ref{uniqueCommonShortest}. Thus
          $f=g=\id_X$. Similarly, $f=g$ is an
          isomorphism of $(X,c_1)$ in the other cases.
    \item $m>n$. This case is similar to the case $m<n$,
          yielding the conclusion that $f=g$ is an
          isomorphism of $(X,c_1)$.
    \item $m=n$. In this case we have either
          $f(B \cup T)= g(B \cup T)=B \cup T$ or
          $f(B \cup T)= g(B \cup T)= L \cup R$. 
          In the former case, $f|_{B \cup T}$ and 
          $g|_{B \cup T}$ are given by one of the
          four possibilities listed above; in the 
          latter case, one of the following holds.
          For $(x,y) \in B \cup T$,
           \begin{itemize}
              \item $f(x,y)=g(x,y)=(y,x)$, or
              \item $f(x,y)=g(x,y)=(y,-x)$, or
              \item $f(x,y)=g(x,y)=(-y,x)$, or
              \item $f(x,y)=g(x,y)=(-y,-x)$.
          \end{itemize}
          An argument like that used above shows that in 
          each of these cases, $f=g$ is an isomorphism
          of $(X,c_1)$.
\end{itemize}
Thus all cases lead to the conclusion that that $f=g$,
hence $A$ is unifying; and that $f \in C(X,c_1)$ 
such that $f(A)=A$ implies $f$ is an isomorphism 
of $(X,c_1)$.
\end{proof}

\subsection{Rectangles in $\Z^2$ with slanted sides and $c_2$}
In this section, we study unifying sets for
digital rectangles with slanted edges in~$\Z^2$,
using the $c_2$ adjacency.
Our assertions are dual to those of
section~\ref{c1RectangleSec} and have proofs with common elements.

\begin{prop}
\label{c2minRectFreeze}
Let $X$ be a digital rectangle in $\Z^2$ with
slanted edges. Let $B \subset X$. Let $B'$ be
the set of endpoints of edges of $X$.
Then $B$ is a freezing set for
$(X,c_2)$ if and only if $B' \subset B$.
Therefore, $B'$ is the only minimal freezing set
for $(X,c_2)$.
\end{prop}

\begin{proof}
By Theorem~\ref{isoInvariant}, there is no loss of generality in assuming
\[B' = \{(0,0), (m,m), (n,-n), (m+n,m-n)\}
  \mbox{ for some } m,n \in \N.
\]
If $B$ is a freezing set, then by 
Proposition~\ref{essentialc2cold90}, $B' \subset B$.
Since $B'$ is a freezing set by 
Theorem~\ref{slantCorners},
it follows that $B'$ is unique as a minimal freezing set.
\end{proof}

\begin{thm}
Let $X$ be the digital rectangle with edges in
the set 
\[ B = \{(0,0), (m,m), (n,-n), (m+n,m-n)\}.
\]
Then $B$ is a unifying set for $(X,c_2)$. Further,
every $f \in C(X,c_2)$ such that $f(B)=B$ is an
isomorphism.
\end{thm}

\begin{proof}
Let $LR$ (lower right) be the edge of $X$ from $(n,-n)$
to $(m+n,m-n)$. Let $UL$ (upper left) be the edge of 
$X$ from $(0,0)$ to $(m,m)$. Let $LL$ (lower left)
be the edge of $X$ from $(0,0)$ to $(n,-n)$. Let
$UR$ (upper right) be the edge of $X$ from $(m,m)$
to $(m+n,m-n)$. For $m<n$, there are distinct isomorphisms
$F_1,F_2,F_3,F_4: S \to S$, where 
\[ S = LR \cup UL \cup LL \cup UR
\]
is the bounding curve
of $X$, where $F_1=\id_S$, $F_2$ reverses the orientations
of $UL$ and $LR$, $F_3$ interchanges $UL$ and $LR$ while
preserving their orientations, and $F_4$ interchanges $UL$ 
and $LR$ and reverses their orientations.

Consider the following cases.
\begin{itemize}
    \item $m<n$. Since $f(B)=g(B)=B$, we have
          that $f(UL)$, $g(UL)$, $f(LR)$, and $g(LR)$
          are $c_2$-paths of length at most $m$ between 
          distinct members of $B$, and since the closest
          distinct members of $B$ are joined by paths
          of length $m$, $f(UL)$, $g(UL)$, $f(LR)$,
          and $g(LR)$ are paths of length $m$. Therefore,
          $f(UL \cup LR)=g(UL \cup LR)=UL \cup LR$.
          Proposition~\ref{uniqueCommonShortest} 
          implies that for all $(x,y)\in UL \cup LR$,
          $f(x,y)=g(x,y) = F_i(x,y)$ for some index~$i$.
          
          Suppose the first case, 
          $f(x,y)=g(x,y)=(x,y)$ for $(x,y)\in UL \cup LR$.
          Each $(x,y) \in X$ lies on the unique shortest
          $c_2$-path (a slanted path) between some
          $d_1 \in UL$ and some $d_2 \in LR$.
          Since $f(d_j)=g(d_j)=d_j$ for $j \in \{1,2\}$,
          we must have $f(x,y)=g(x,y)=(x,y)$ by
          Proposition~\ref{uniqueCommonShortest}. Thus
          $f=g=\id_X$. Similarly, $f=g$ is an
          isomorphism of $(X,c_2)$ in the other cases.
    \item $m>n$. This case is similar to the case $m<n$,
          and we similarly conclude that $f=g$ is an
          isomorphism of $(X,c_2)$.
    \item $m=n$. Here, in addition to the isomorphisms
          $F_1,F_2,F_3,F_4$ discussed above, we also
          have isomorphisms $R_1,R_2,R_3,R_4$ of $(X,c_2)$
          that rotate the edges of $X$ by $90^{\circ}$
          ($\pi /2$ radians) either clockwise or
          counterclockwise, either preserving or
          reversing the orientations of both $UL$
          and $LR$. An argument like that used above 
          shows that in each of these cases, $f=g$ 
          is an isomorphism of $(X,c_2)$.
\end{itemize}
Thus all cases lead to the conclusion that $f=g$,
hence $B$ is unifying; and that $f \in C(X,c_2)$ 
such that $f(B)=B$ implies $f$ is an isomorphism
of $(X,c_2)$.
\end{proof}

\subsection{Generalized normal product}
In this section, we consider unifying sets for
Cartesian products of digital images using the
normal product adjacency.

We have the following generalization of 
the normal product adjacency~\cite{Berge}
for the Cartesian product of two graphs.

\begin{definition}
{\rm \cite{Sabidussi,BxNormal}}
Let $u,v \in \N$, $1 \le u \le v$.
Let $(X_i,\kappa_i)$ be digital images,
$i \in \{1,\ldots,v)$. Let $x_i,y_i \in X_i$,
$x = (x_1, \ldots, x_v)$, $y=(y_1,\ldots,y_v)$.
Then $x \adj y$ in the 
{\em generalized normal product adjacency}
$NP_u(\kappa_1, \ldots, \kappa_v)$ if
for at least 1 and at most $u$ indices $i$,
$x_i \adj_{\kappa_i} y_i$ and for all other
indices $j$, $x_j=y_j$.
\end{definition}

\begin{remk}
For $u=v=2$, the generalized normal product adjacency 
coincides with the {\em normal product adjacency}.
Sabidussi~\cite{Sabidussi} uses {\em strong} for
what we call the generalized normal product adjacency;
we prefer the latter name, as ``strong" also appears
in the literature for what we call the normal product
adjacency.
\end{remk}

The following generalizes a result
in~\cite{BK12,Boxer16}.

\begin{thm}
\label{prod-cont}
{\rm \cite{BxNormal}}
Let $f_i: (X_i, \kappa_i) \to (Y_i, \lambda_i)$, $1 \le i \leq v$.
Then the product map
\[ f=\Pi_{i=1}^v f_i: (\Pi_{i=1}^v X_i, NP_v(\kappa_1, \ldots, \kappa_v)) \to (\Pi_{i=1}^v Y_i, NP_v(\lambda_1, \ldots, \lambda_v)) \]
given by $f(x_1, \ldots, x_v) = (f_1(x_1), \ldots, f_v(x_v))$
is continuous if and only if each $f_i$ is continuous.
\end{thm}

\begin{thm}
Let $\emptyset \neq A_i \subset X_i$,
where $(X_i, \kappa_i)$ is a digital image,
$1 \le i \le v \in \N$. Let
$A = \Pi_{i=1}^v A_i$, $X = \Pi_{i=1}^v X_i$.
If $A$ is a unifying set for 
$(X,NP_v(\kappa_1, \ldots, \kappa_v))$ then
for each $i$, $A_i$ is a unifying set for 
$(X_i, \kappa_i)$.
\end{thm}

\begin{proof}
Suppose $A$ is a unifying set for 
$(X,NP_v(\kappa_1, \ldots, \kappa_v))$. For all~$i$, let
$f_i,g_i \in C(X_i,\kappa_i)$ be such that
$f_i(A_i)=g_i(A_i)=A_i$ and $f_i|_{A_i} = g_i|_{A_i}$.
Then by Theorem~\ref{prod-cont},
$f = f_1 \times \cdots \times f_v$ and
$g = g_1 \times \cdots \times g_v$ are members of
$C(X,NP_v(\kappa_1, \ldots, \kappa_v))$. 
Further, given $a = (a_1,\ldots,a_v) \in A$,
there exist $a_i' \in A_i$ such that
$f_i(a_i') = g_i(a_i')=a_i$, and therefore we have
$f(A) = g(A) = A$ and $f|_A=g|_A$. Since $A$ is
unifying, we have $f=g$, and therefore $f_i=g_i$
for all~$i$. Thus $A_i$ is unifying.
\end{proof}

\section{Shy maps that are retractions}
Shy maps in digital topology were introduced in~\cite{Bx05} 
and studied further
in~\cite{Bx14,BxSt16,Boxer16,BxNormal,BxAlternate}. A version
of shy maps for topological spaces was introduced
in~\cite{BxTopShy}.

\begin{definition}
\label{shy-def}
\cite{Bx05}
Let $f: (X,\kappa) \rightarrow (Y,\lambda)$ be a
continuous function of digital images. We say $f$ is
{\em shy} if
\begin{itemize}
\item for each $y \in f(X)$, $f^{-1}(y)$ is connected, and
\item for every $y_0,y_1 \in f(X)$ such that $y_0$ and $y_1$ are
      adjacent, $f^{-1}(\{y_0,y_1\})$ is connected. 
\end{itemize}
\end{definition}

\begin{thm}
\label{shy-thm}
{\rm ~\cite{Bx14}}
Let $f: (X,\kappa) \rightarrow (Y,\lambda)$ be a 
continuous function between digital images.
Then $f$ is shy if and only if
$f^{-1}: Y \multimap X$ is a connectivity preserving multivalued
function (i.e., given $\lambda$-connected $A \subset f(Y)$,
$f^{-1}(A)$ is $\kappa$-connected).
\end{thm}

We say a point $p$ of a connected graph $G=(X,\kappa)$ is an
{\em articulation point} of $G$ if $(X \setminus \{p\},\kappa)$
is not connected.

\begin{thm}
\label{shyRetraction}
Let $(X,\kappa)$ be a digital image.
Let $\emptyset \neq R \subset X$. Let 
$\emptyset \neq A \subset R$ such that
for each $\kappa$-component $K$ of $X \setminus R$
there exists $p \in A$ such that $p$ is an articulation point for $(K \cup R,\kappa)$. Then there
is a unique function $r: X \to R$ that is
a shy $\kappa$-retraction.
\end{thm}

\begin{proof}
For $x \in X \setminus R$, let $p_x \in A$ be the
articulation point for the union of $R$ and the
$\kappa$-component $K_x$ of $X \setminus R$ 
containing $x$.
Let $r: X \to X$ be the function
\[ r(x)= \left \{ \begin{array}{ll}
    x & \mbox{if } x \in R; \\
    p_x & \mbox{if } x \in X \setminus R.
\end{array} \right .
\]
Clearly, $r(X)=R$ and $r|_R = \id_R$. It is easily
seen that $r^{-1}(p_x) \setminus \{p_x\} = K_x$ 
is a $\kappa$-component of
$X \setminus R$, and $r^{-1}(y)=\{y\}$ for
$y \in R \setminus T$. It follows that $r \in C(X,\kappa)$,
that $r$ is a retraction of $X$ to $R$, and $r$ is shy.

Suppose $f \in C(X,\kappa)$ is a shy retraction of 
$X$ to $R$. If there exists $x_0 \in X \setminus R$
such that $x_1 = f(x_0) \neq p_{x_0}$, then $p_{x_0}$ separates 
the points $x_0,x_1 \in f^{-1}(x_1)$,
contrary to the assumption that $f$ is shy.
The uniqueness of $r$ as a shy retraction follows.
\end{proof}

\begin{cor}
Let $(X,\kappa)$ be a digital image that is a tree.
Let $(R,\kappa)$ be a nonempty subtree of $(X,\kappa)$.
Then there is a unique function $r: X \to R$ that is
a shy $\kappa$-retraction.
\end{cor}

\begin{proof}
It is trivial that if $R = X$, we can take $r = \id_X$.
Otherwise, we can take
\[A = \{x \in R ~|~ x 
     \mbox{ is a leaf of $R$ and not a leaf of } X\}.
\]
The assertion follows from Theorem~\ref{shyRetraction}.
\end{proof}

For topological spaces, we have the following.

\begin{definition}
\label{shy-def}
{\rm \cite{BxTopShy}}
Let $X$ and $Y$ be topological spaces and let $f: X \to Y$. Then
$f$ is {\em shy} if $f$ is continuous and for every path-connected $Y' \subset f(X)$, 
$f^{-1}(Y')$ is a path-connected subset of $X$. $\Box$
\end{definition}

By using an argument similar to the proof of
Theorem~\ref{shyRetraction}, we get the following.

\begin{thm}
Let $X$ be a topological space.
Let $\emptyset \neq A \subset R \subset X$
such that each $p \in A$ separates $R$ and a
component of $X \setminus R$. Then there
is a unique continuous function $r: X \to R$ that is
a shy retraction.
\end{thm}

\section{Approximate fixed points}
Suppose $A \subset X$ and $A$ is a $\kappa$-freezing set
for $X$. By definition, if $f \in C(X,\kappa)$ and 
$A \subset \Fix(f)$, then $f = \id_X$, i.e., $X = \Fix(f)$.
If we weaken the hypothesis so that
instead of assuming $A \subset \Fix(f)$ we assume every
point of $A$ is an approximate fixed point of $f$, might we
reach the weaker conclusion that every point of $X$
is an approximate fixed point of $f$? The answer is not generally
affirmative; we give a counterexample below. We also examine 
basic examples for which an affirmative answer is shown.

\subsection{Wedge of cycles}
In this section, we show that a wedge of cycles $X$ can
support a freezing set $A$ and a continuous self-map $f$
such that every point of $A$ is an approximate fixed point of $f$,
but not every point of $X$ is an approximate fixed point of $f$.

\begin{thm}
\label{wedgeSCCthm}
{\rm \cite{BxFpSets}}
Let $C_m$ and $C_n$ be cycles, with $m>4$, $n>4$, where
$C_m = \{x_i\}_{i=0}^{m-1}$, $C_n=\{x_i'\}_{i=0}^{n-1}$, with
the members of $C_m$ and $C_n$ indexed circularly.
Let $x_0=x_0'$ be the wedge point of $X=C_m \vee C_n$. Let
$x_i,x_j \in C_m$ and $x_k',x_p' \in C_n$ be such that
$C_m$ is the union of unique shorter paths determined by
$x_i, x_j, x_0$ and $C_n$ is the union of unique 
shorter paths determined by $x_k', x_p', x_0'$. Then
$A=\{x_i, x_j, x_k', x_p'\}$ is a freezing set for $X$.
\end{thm}

\begin{figure}
    \centering
    \includegraphics[height=1.5in]{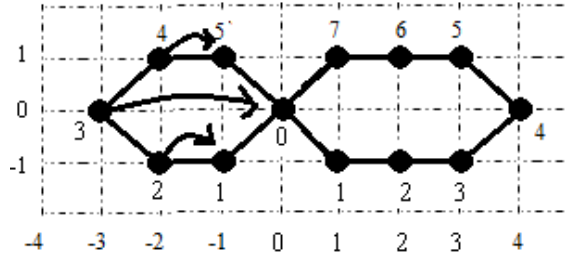}
    \caption{\cite{BxConvex} The map $f$ of 
    Example~\ref{wedgeExl}. Points are labeled
    by their indices as in the Example. The cycle with points
    $p=(x,y)$ for $x \le 0$ represents $C_6$, for which
    $\{x_0,x_2,x_4\}$ is a $c_2$-freezing set; the cycle with points
    $p=(x,y)$ for $x \ge 0$ represents $C_m$ (here, $m=8$, and
    $\{x_0',x_3',x_6'\}$ is a $c_2$-freezing set for $C_8$, so
    $A=\{x_2,x_4,x_3',x_6'\}$ is a $c_2$-freezing set for
    $C_6 \vee C_8$). Arrows connect $p$ and $f(p)$ for points
    $p \not \in \Fix(f)$. Each point of $A$ is a $c_2$-approximate
    fixed point of~$f$.
    }
    \label{fig:wedgeCyclesWMap}
\end{figure}

\begin{exl}
\label{wedgeExl}
Let $X = C_6 \vee C_m$, where
$C_6 = \{x_i\}_{i=0}^{5}$ and $C_n=\{x_i\}_{i=0}^{n-1}$ are
$c_2$-simple closed curves in $\Z^2$, with
the members of $C_6$ and $C_n$ indexed circularly.
By Theorems~\ref{3ptsForCycles} and~\ref{wedgeSCCthm}, if
$k$ and $p$ are chosen so
that $\{x_0', x_k', x_p'\}$ is a freezing set for $C_n$,
then we can take $A = \{x_2, x_4, x_k', x_p'\}$ to be a freezing
set for $(X,c_2)$. Now take $f: X \to X$ to be the function
\[ f(x) = \left \{ \begin{array}{ll}
    x_0 & \mbox{if } x = x_3; \\
    x_1 & \mbox{if } x = x_2; \\
    x_5 & \mbox{if } x = x_4; \\
    x   & \mbox{otherwise.}
\end{array} \right .
\]
See Figure~\ref{fig:wedgeCyclesWMap}. One sees easily that $f \in C(X,c_2)$, that every member of
$A$ is a $c_2$-approximate fixed point of $f$, but $x_3$ is not
a $c_2$-approximate fixed point of~$f$.
\end{exl}

\subsection{Disks in $(\Z^2,c_1)$}
\label{c1Disks}
\begin{lem}
\label{segmentLemma}
Let $q_0,q_1 \in X \subset \Z^2$. Suppose there is a horizontal
or vertical $c_1$-path $P$ in $X$ from $q_0$ to $q_1$. Let
$f: P \to X$ be $c_1$-continuous, such that $q_0$ and
$q_1$ are $c_1$-approximate fixed points of $f$. Then every
member of $P$ is a $c_1$-approximate fixed point of $f$.
\end{lem}

\begin{proof}
Without loss of generality, $P$ is horizontal, $q_0 = (0,0)$,
and $q_1 = (n,0)$ for some $n \in \N$. Suppose there exists
$q = (x,0) \in P$ such that $q$ is not a $c_1$-approximate fixed
point of $f$. Then $|x - p_1(f(q))|>1$; or $|p_2(f(q))|>1$;
or $|x - p_1(f(q))|=1$ and $|p_2(f(q))|=1$.

If $|x - p_1(f(q))|>1$ then either $p_1(f(q))>x+1$ or
$p_1(f(q)) < x-1$.
\begin{itemize}
    \item Suppose $p_1(f(q)) > x+1$. Then by Lemma~\ref{c1pulling}
          we would have $p_1(f(q_0)) > 1$, contrary to the
          assumption that $q_0 = (0,0)$ is an approximate
          fixed point.
    \item If $p_1(f(q)) < x-1$, then by Lemma~\ref{c1pulling}
          we would have $p_1(f(q_1)) < n-1$, contrary to
          the assumption that $q_1 = (n,0)$ is an approximate
          fixed point.
\end{itemize}

Suppose $|p_2(f(q))|>1$. Without loss of generality,
$p_2(f(q))>1$, as the case $p_2(f(q))<1$ can be handled
similarly. Since $c_1$-adjacent points differ in only one
coordinate and the $q_i$ as approximate fixed points 
implies $|p_2(f(q_i)| \le 1$, $i \in \{0,1\}$, there are
at least 4 indices~$j$ for which $p_2(f(x_j)) \neq p_2(f(x_{j+1}))$
and therefore at most $n-4$ indices~$j$ for which 
$p_1(f(x_j)) \neq p_1(f(x_{j+1}))$. This is a contradiction,
since $x_0$ and $x_1$ being approximate fixed points 
implies $p_1(f(x_0)) \le 1$ and $p_1(f(x_1) \ge n-1$, so
at least $n-2$ indices~$j$ would satisfy
$p_1(f(x_j)) \neq p_1(f(x_{j+1}))$.

Suppose $|x - p_1(f(q))|=1$ and $|p_2(f(q))|=1$. Without loss
of generality, $p_1(f(q))=x+1$ and $p_2(f(q))=y+1$. By the
$c_1$-continuity of $f$ and Lemma~\ref{c1pulling} it follows
that $p_1(f(q_0))\ge 1$ and $p_2(f(q_0))\ge 1$, contrary to
the assumption that $q_0$ is a $c_1$-approximate fixed point
of~$f$.

Thus every case yields a contradiction brought about by
assuming there is a point of $P$ that is not an
approximate fixed point of $f$. The assertion follows.
\end{proof}

\begin{thm}
\label{AFPc1}
Let $V_i \subset X \subset \Z^2$, $i \in \{1,\ldots,n\}$ 
where each $V_i$ is a thick convex disk. 
Let $X' = \bigcup_{i=1}^n V_i$.
Let $C_i$ be a bounding curve of $V_i$.
Let $A_{1,i}$ be the set of endpoints of maximal
axis-parallel segments of $C_i$. Let $A_{2,i}$ be the 
union of maximal slanted segments of $C_i$. 
\begin{enumerate}
    \item $A = (X \setminus X') \cup \bigcup_{i=1}^n (A_{1,i} \cup A_{2,i})$ is a freezing set for $(X,c_1)$.
    \item Suppose $f \in C(X,c_1)$ such that every point of $A$ is
a $c_1$-approximate fixed point of $f$. Then every point of $X$
is a $c_1$-approximate fixed point of $f$.
\end{enumerate}
\end{thm}

\begin{proof}
Assertion 1) is Theorem~\ref{corners}. To prove 
assertion 2), we argue as follows.

Let $S$ be a maximal digital segment of a bounding curve $C_i$
for $V_i$. If $S$ is horizontal or vertical, then
by Lemma~\ref{segmentLemma}, every point of $S$ is
a $c_1$-approximate fixed point of $f$.
If $S$ is slanted, then $S \subset A$, so every point of $S$ is
a $c_1$-approximate fixed point of $f$.
Thus each point of $C_i$, is
a $c_1$-approximate fixed point of $f$.

For $x \in X \setminus A$, there is a horizontal segment $P$
containing $x$ such that the endpoints of $P$ belong to 
$\bigcup_{i=1}^n C_i$,
and therefore are approximate fixed points of $f$. By
Lemma~\ref{segmentLemma}, every point of $P$ is 
a $c_1$-approximate fixed point of $f$. Thus, every point of
$X$ is a $c_1$-approximate fixed point of $f$.
\end{proof}

\begin{remk}
Theorems~\ref{corners} and~\ref{AFPc1} simplify when $X'=X$, 
in which case $A=\bigcup_{i=1}^n (A_{1,i} \cup A_{2,i})$. They
might be applied in this case when $i \neq j$ implies
$V_i \cap V_j$ is empty, a single point, or a common edge
of $V_i$ and $V_j$.
\end{remk}

\subsection{Disks in $(\Z^2,c_2)$}
We show in this section that disks in $(\Z^2,c_2)$ yield
results similar to those shown in section~\ref{c1Disks} for
the $c_1$ adjacency. 

\begin{lem}
\label{segmentLemmaC2}
Let $q_0,q_1 \in X \subset \Z^2$. Suppose there is a slanted
$c_2$-path $P$ in $X$ from $q_0$ to $q_1$. Let
$f: P \to X$ be $c_2$-continuous, such that $q_0$ and
$q_1$ are $c_2$-approximate fixed points of $f$. Then every
member of $P$ is a $c_2$-approximate fixed point of $f$.
\end{lem}

\begin{proof}
Without loss of generality, the slope of $P$ is 1.
Without loss of generality, $q_0 = (0,0)$ and 
$q_1 = (n,n)$ for $n=length(P)$.
Suppose there exists $p \in P$ that is not a $c_2$-approximate 
fixed point of $f$. Then $|p_1(f(p)) - p_1(p)|>1$ or
$|p_2(f(p)) - p_2(p)|>1$.
\begin{itemize}
    \item If $|p_1(f(p)) - p_1(p)|>1$ then either
          $p_1(f(p)) - p_1(p)>1$ or $p_1(p) - p_1(f(p))>1$.
          \begin{itemize}
              \item If  $p_1(f(p)) - p_1(p)>1$ then by
                    Lemma~\ref{c1pulling},
                    $1 < p_1(f(q_0)) - p_1(q_0) = p_1(f(q_0))$,
                    contrary to the assumption that $q_0$ is an
                    approximate fixed point.
              \item If $p_1(p) - p_1(f(p))>1$, then by
                    Lemma~\ref{c1pulling},
                    $1 < p_1(q_1) - p_1(f(q_1)) = n - p_1(f(q_1))$,
                    or $p_1(f(q_1)) < n-1$, 
                    contrary to the assumption that $q_1$ is an
                    approximate fixed point.
          \end{itemize}
    \item If $|p_2(f(p)) - p_2(p)|>1$ then, similarly, we obtain
          contradictions.
\end{itemize} 
Since all cases yield contradictions, the hypothesis of a
$p \in P$ that is not a $c_2$-approximate fixed point of $f$
must be false. This completes the proof.
\end{proof}

The following is a dual to Theorem~\ref{AFPc1}.

\begin{thm}
\label{AFPc2}
Let $V_i \subset X \subset \Z^2$, $i \in \{1,\ldots,n\}$ 
where each $V_i$ is a thick convex disk. 
Let $X' = \bigcup_{i=1}^n V_i$.
Let $C_i$ be a bounding curve of $V_i$.
Let $B_{1,i}$ be the union of maximal horizontal 
and maximal vertical segments of $C_i$.
Let $B_{2,i}$ be the set of endpoints of maximal slanted
segments of $C_i$. 
\begin{enumerate}
    \item $B = (X \setminus X') \cup \bigcup_{i=1}^n (B_{1,i} \cup B_{2,i})$ is a freezing set for $(X,c_2)$.
    \item Suppose $f \in C(X,c_2)$ such that every point of $B$ is
a $c_2$-approximate fixed point of $f$. Then every point of $X$
is a $c_2$-approximate fixed point of $f$.
\end{enumerate}
\end{thm}

\begin{figure}
    \centering
    \includegraphics[height=1.5in]{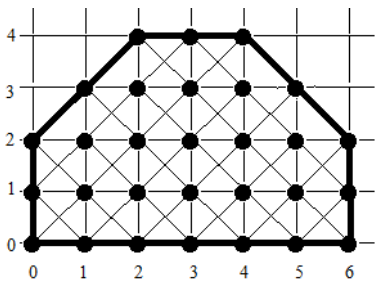}
    \caption{The point (3,2) in the digital image shown above
    does not lie on a slanted segment that joins 2 points of the
    boundary curve shown darkly.
    }
    \label{fig:slantToBdNotMeetBdPt}
\end{figure}

\begin{proof}
Assertion 1) is Theorem~\ref{slantCorners}. To
prove assertion 2), we argue as follows.

By Lemma~\ref{segmentLemmaC2}, every slanted segment of $C_i$ is
made up entirely of $c_2$-approximate fixed points of $f$.
From Theorem~\ref{convDiskThmC2Actual} it follows that $C_i$ is
made up entirely of $c_2$-approximate fixed points of $f$.
Therefore, every point of $B$ is a $c_2$-approximate fixed point
of $f$.

Lemma~\ref{segmentLemmaC2} lets us conclude that
if $x \in X$ such that $x$ lies on a slanted segment $P$ that
connects two points of $B$, then $x$ is a $c_2$-approximate
fixed point of $f$.

This leaves us to consider points $p=(x_0,y_0) \in X$ such that $p$ 
does not lie either on an axis-parallel segment of $B$ or
on a slanted segment $P$ that
connects two points of $B$. Such a point must be in the
interior of $X$ and therefore is $c_2$-adjacent to its
4 $c_1$-neighbors $q_1=(x_0-1,y_0)$, $q_2=(x_0+1,y_0)$, 
$q_3=(x_0, y_0-1)$, and $q_4=(x_0,y_0+1)$, each of which
lies on a slanted segment joining members of $S$
(see Figure~\ref{fig:slantToBdNotMeetBdPt}). Therefore,
by Lemma~\ref{segmentLemmaC2}, $q_1$, $q_2$, $q_3$, and $q_4$ are
approximate fixed points of $f$.

Suppose $p$ is not a $c_2$-approximate fixed point of $f$.
Then either $|p_1(f(p)) - x_0| > 1$ or $|p_2(f(p)) - y_0| > 1$.
\begin{itemize}
    \item Suppose $|p_1(f(p)) - x_0| > 1$. Then either
          $p_1(f(p)) - x_0 > 1$ or $x_0 - p_1(f(p)) > 1$.
          \begin{itemize}
              \item Suppose $p_1(f(p)) - x_0 > 1$. Then by the
                    continuity of $f$ and Lemma~\ref{c1pulling},
                    $p_1(q_1) - p_1(f(q_1)) > 1$, contrary to
                    $q_1$ being an approximate fixed point of $f$.
              \item Suppose $x_0 - p_1(f(p)) > 1$. Then by the
                    continuity of $f$ and Lemma~\ref{c1pulling},
                    $p_1(q_2) - p_1(f(q_2)) > 1$, contrary to
                    $q_2$ being an approximate fixed point of $f$.
          \end{itemize}
    \item Similarly, we obtain a contradiction if
          $|p_2(f(p)) - y_0| > 1$.
\end{itemize}
Since all cases yield a contradiction when we assume
$p$ is not a $c_2$-approximate fixed point of $f$, this
hypothesis must be incorrect. The assertion follows.
\end{proof}

\begin{remk}
Theorems~\ref{slantCorners} and~\ref{AFPc2} simplify when $X'=X$, 
in which case $B=\bigcup_{i=1}^n (B_{1,i} \cup B_{2,i})$. They
might be applied in this case when $i \neq j$ implies
$V_i \cap V_j$ is empty, a single point, or a common edge
of $V_i$ and $V_j$.
\end{remk}

\subsection{Trees}
\begin{thm}
{\rm \cite{BxFpSets}}
\label{tree}
Let $(X,\kappa)$ be a digital image such that
the graph $G=(X,\kappa)$ is a finite tree with $\#X > 1$.
Let $E$ be the set of vertices of $G$ that have degree 1. Then
$E$ is a minimal freezing set for $G$.
\end{thm}

\begin{lem}
\label{keepSegmentInside}
Let $(X,\kappa)$ be a digital image such that
the graph $G=(X,\kappa)$ is a finite tree. Let $f \in C(X,\kappa)$.
Let $a,b \in X$ be such that $a$ and $b$ are
$\kappa$-approximate fixed points of $f$. Let $P$ be the
unique shortest path in $G$ from $a$ to $b$. Then
$f(P) \subset P$ and every point of $P$ is 
 a $\kappa$-approximate fixed point of $f$.
\end{lem}

\begin{proof}
Let $P = \{x_i\}_{i=0}^n$ such that $x_0=a$, $x_n=b$, and
$x_i \adj_{\kappa} x_j$ if and only if $|i-j|=1$.
\begin{itemize}
    \item Suppose $f(a)=a$. Let us show that
          \begin{equation} 
          \label{penultimate}
          f(x_{n-1}) \in \{x_{n-1}, b\} \subset P.
          \end{equation}
          We know that
          $f(b) \in N^*(X,b,\kappa)$. If~(\ref{penultimate}) is
          false, then $f(P)=P \cup \{f(b)\}$ is the unique shortest
          path in $G$ from $a$ to $f(b)$. But $P \cup \{f(b)\}$
          has length $n+1$, and $\#P = n+1$ implies
          $\#f(P) \le n$, so we have a contradiction brought about
          by negating~(\ref{penultimate}). Thus~(\ref{penultimate})
          is established.
          
          Since $G$ is acyclic, we must have $f(P) \subset P$. Now
          suppose for some $k$ that $x_k$ is not an approximate
          fixed point of $f$. Then $f(x_k) = x_m$ for some $m$
          such that $|k-m| > 1$. Without loss of generality,
          $m-k>1$. Then by continuity and since $G$ is acyclic,
          $f(x_k)$ must ``pull"~\cite{hmps} $f(a)=f(x_0)$ so that
          $f(a) = x_t$ for some $t > 1$, contrary to
          $a \in \Fix(f)$. The contradiction establishes
          that each point of $P$ must be an approximate
          fixed point of $f$.
    \item Suppose $f(a) \not \in P$. Recall we are assuming
          $f(b) \in N^*(X,b,\kappa)$, so $f(b) \in \{x_{n-1},b\}$
          or $f(b) \not \in P$. We claim $f(b)=x_{n-1}$.
          For otherwise, 
          $f(P) = \{f(a) \neq x_0, a=x_0, x_1, \ldots, x_n=b, f(b)\}$
          where $f(b)$ may be equal to $b$, so 
          $\#f(P) \in \{n+2, n+3\}$ while $\#P=n+1$,
          a contradiction.
          Therefore, $f(b)=x_{n-1}$. By the acyclicity of $G$,
          it follows that 
          $f(P) = \{f(a)\} \cup \{x_i\}_{i=0}^m$, where
          $m \in \{n-1,n\}$. As in the case
          $f(a)=a$, it follows that every point of $P$ is an
          approximate fixed point of $f$.
    \item Suppose $f(a) \in P \setminus \{a\}$. Since $a$ is an
          approximate fixed point of $f$, it follows that 
          $f(a)=x_1$. It follows as in the case $f(a) \not \in P$
          that every point of $P$ is an
          approximate fixed point of $f$.
\end{itemize}
This establishes the assertion.
\end{proof}

\begin{thm}
\label{treeFreezeEssential}
Let $(X,\kappa)$ be a digital image such that
the graph $G=(X,\kappa)$ is a finite tree with $\#X > 1$.
Let $E$ be the set of vertices of $G$ that have degree 1. Then
given a freezing set $A$ for $G$, we have $E \subset A$.
\end{thm}

\begin{proof}
Since $\#X > 1$, we can choose $x_0 \in E \setminus \{a\}$ to
be a root for $G$. Then the function $f: X \to X$ given by
\[ f(x) = \left \{ \begin{array}{ll}
    x & \mbox{for } x \neq a; \\
    parent(a) & \mbox{for } x = a, 
\end{array}
\right .
\]
is easily seen to be a member of $C(X,\kappa)$ such that
$f|_{E \setminus \{a\}} = \id_{E \setminus \{a\}}$ and
$f(a) \neq a$. Thus $E \setminus \{a\}$ cannot be a
freezing set for $G$. The assertion follows.
\end{proof}

\begin{thm}
Let $(X,\kappa)$ be a digital image such that
the graph $G=(X,\kappa)$ is a finite tree with $\#X > 1$.
Let $A$ be a freezing set for $G$. Suppose
$f \in C(X,\kappa)$ is such that for each $a \in A$,
$a$ is an approximate fixed point of $f$. Then for all $x \in X$,
$x$ is an approximate fixed point of $f$.
\end{thm}

\begin{proof}
Let $E$ be the set of vertices of $G$ that have degree 1.
By Theorem~\ref{treeFreezeEssential}, $E \subset A$, and by
Theorem~\ref{tree}, $E$ is a freezing set. Therefore,
there is no loss of generality in assuming $A = E$.

Let $f \in C(X,\kappa)$ such that
for each $e \in E$, $e$ is an approximate fixed point of $f$.
We can choose $x_0 \in E$ as a root of $X$.
Since $x \in X$ implies $x$ is on the unique shortest path 
in~$G$ from $x_0$ to some $e \in E$, it follows from
Lemma~\ref{keepSegmentInside} that $x$ is an approximate 
fixed point of $f$.
\end{proof}

\subsection{Cycles}

\begin{thm}
Let $(C_n, \kappa)$ be a digital cycle of $n$ distinct points, 
$n \in \N$, $n \ge 3$, with $C_n = \{x_i\}_{i=0}^{n-1}$,
such that $x_i \adj_{\kappa} x_j$ if and only if 
$j = (i \pm 1) \mod n$. Let
$A = \{x_u, x_v, x_w\}$ be a set of distinct members of $C_n$
such that $C_n$ is a union of unique shorter paths determined 
by these points. Let $f \in C(C_n, \kappa)$ be such that
every member of $A$ is an approximate fixed point of $f$.
Then every member of $C_n$ is an approximate fixed point of $f$,
and $f$ is an isomorphism.
\end{thm}

\begin{proof} Note by Theorem~\ref{3ptsForCycles}, $A$
is a minimal freezing set for $(C_n, \kappa)$.

First, we show that f must be a surjection. Without loss of
generality, $0 \le u < v < w < n$. Suppose $B$ is the
unique shorter path in $C_n$ from $x_u$ to $x_v$. Since we
must have $\#f(B) \le \#B$ and $x_u$ and $x_v$ are approximate
fixed points, we must have $f(x_u) \in \{x_{u-1}, x_u, x_{u+1}\}$
and $f(x_v) \in \{x_{v-1}, x_v, x_{v+1}\}$.

 Suppose $f(x_u) = x_u$. We must have 
          \[ \#f(B) \le \#B = v-u+1 \le n/2,
          \] so 
          $f(x_v) \in \{x_{v-1},x_v\}$. If $f(x_v)=x_{v-1}$,
          then we must have $f(x_w)=x_{w-1}$, hence
          (proceeding with increasing indices, $\mod n$,)
          $f(x_{u-1})=x_{u-2}$, so $f$ would be discontinuous
          at the adjacent pair $x_{u-1}$ and $x_u$. Thus we would
          have $f(x_v) = x_v$ and $f(x_w)=x_w$. Thus $f|_A = \id_A$.
Since $A$ is freezing, it follows that $f = \id_X$.

If $f(x_u)=x_{u-1}$ or $f(x_u)=x_{u+1}$, we can apply a
rotation $r(x_i) = x_{(i-1) \mod n}$ (respectively,
$r(x_i) = x_{(i+1) \mod n}$, which is an isomorphism.
Then by the above, $r \circ f = \id_X$ is an isomorphism, so
\[ f = r^{-1} \circ r \circ f =  r^{-1} \circ \id_X = r^{-1}
\]
is an isomorphism, with each member of $A$ an approximate fixed
point.

Thus, in all cases, each member of $A$ is an approximate fixed
point of $f$, which must be an isomorphism.
\end{proof}

\section{Further remarks}
When a member of $C(X,\kappa)$ has restricted behavior on
a subset $A$ of $X$, the restriction may have a powerful
effect on the behavior of $f|_{X \setminus A}$. We have
examined instances of this phenomenon with respect to
freezing and cold sets, retractions, and shy maps, on a
variety of basic digital images.

\end{document}